\documentclass{amsart}

\theoremstyle{definition}

\theoremstyle{remark}

\numberwithin{equation}{section}

\begin{document}

\title{Landen Transformations and the Integration of Rational Functions}

\author{George Boros}
\address{Department of Mathematics, Xavier University, New 
Orleans, Louisiana 70125}
\email{gboros@xula.math.edu}

\author{Victor H. Moll}
\address{Department of Mathematics,
Tulane University, New Orleans, LA 70118}
\email{vhm@math.tulane.edu}
\thanks{The second author was supported in part 
by NSF Grant DMS-0070567.}

\subjclass{Primary 33}

\date{\today}


\begin{abstract}
We present a rational version of the classical Landen transformation for 
elliptic integrals.  This is employed to obtain explicit closed-form 
expressions for a large class of integrals of even rational functions 
and to develop an algorithm for numerical integration of these functions. 
\end{abstract}

\maketitle


\newtheorem{Definition}{\bf Definition}[section]
\newtheorem{Thm}[Definition]{\bf Theorem} 
\newtheorem{Lem}[Definition]{\bf Lemma} 
\newtheorem{Cor}[Definition]{\bf Corollary} 
\newtheorem{Prop}[Definition]{\bf Proposition}

\newcommand{\nn}{\nonumber}
\newcommand{\ba}{\begin{eqnarray}}
\newcommand{\bh}{B \left( r - \tfrac{1}{2}, \tfrac{1}{2} \right)}
\newcommand{\ea}{\end{eqnarray}}
\newcommand{\al}{\beta}
\newcommand{\be}{\beta}
\newcommand{\dz}{\frac{d}{dz}}
\newcommand{\exta}{c}
\newcommand{\E}{{\mathfrak{E}}}
\newcommand{\rat}{{\mathfrak{R}}}
\newcommand{\F}{{\mathfrak{F}}}
\newcommand{\Ro}{{\mathfrak{R}}}
\newcommand{\ift}{\int_{0}^{\infty}}
\newcommand{\no}{\noindent}
\newcommand{\ga}{a}
\newcommand{\x}{a}
\newcommand{\m}{m}
\newcommand{\X}{{\mathbb{X}}}
\newcommand{\Q}{{\mathbb{Q}}}
\newcommand{\R}{{\mathbb{R}}}
\newcommand{\Y}{{\mathbb{Y}}}
\newcommand{\qq}{x^4 + 2ax^2 + 1}
\newcommand{\qqt}{t^4 + 2at^2 + 1}
\newcommand{\qqe}{x^{8} + a_{2}x^{6} + 2a_{1}x^{4} + a_{2}x^{2} + 1}
\newcommand{\bb}{bx^4 + 2ax^2 + 1}
\newcommand{\rplus}{{\mathbb{R}}_{+}}

\section{Introduction} \label{sec-intro} 
\setcounter{equation}{0}
 
We consider the space of
even rational functions of degree $2p$
\ba
{\E}_{2p} & := & \left\{ R(z) = \frac{P(z)}{Q(z)} \; {\Big{|}} \; 
P(z) := \sum_{k=0}^{p-1} b_{k} z^{2(p-1-k)} \text{ and }
Q(z) := \sum_{k=0}^{p} a_{k} z^{2(p-k)} \right\} \nn
\ea
\no
with positive real coefficients $a_{k}, \, b_{k} \in \rplus$ 
normalized by the condition $a_{0} = a_{p} = 1$,  the space
\ba
{\E}_{\infty}  & := & \bigcup_{p=1}^{\infty} {\E}_{2p} \nn
\ea
\no
of normalized even rational functions, and the set 
of $2p-1$ parameters 
$${\mathfrak{P}}_{2p} := \{ a_{1}, \cdots, a_{p-1}; b_{0}, \cdots, 
b_{p-1} \}.$$ 
\no
We describe an algorithm to determine, as a 
function of the parameter set 
${\mathfrak{P}}_{2p}$,  a {\em closed-form} expression of the 
integral 
\ba
I & := & \ift R(z) \, dz \label{integral1}
\ea
\no
for a large class of functions $R \in {\E}_{\infty}$. The function $R$ is 
called {\em symmetric} if its denominator $Q$ satisfies $Q(1/z) = z^{-2p}Q(z)$.
This is equivalent to its coefficients being palindromic, i.e.
$a_{j} = a_{p-j}$ for $1 \leq j \leq p$.  

The class of symmetric functions plays a crucial role in this 
algorithm.  Define
\ba
{\E}_{2p}^{s} & := & \{ R \in {\E}_{2p} \; {\Big{|}} \;  
den(R) \text{ is symmetric} \}
\nn
\ea 
\no
(where $den(R)$ denotes the denominator of $R$), 
the class of rational functions with symmetric denominators of degree $2p$, and
\ba
{\E}_{\infty}^{s}  & := & \bigcup_{p=1}^{\infty} {\E}_{2p}^{s}. \nn
\ea
\no
For $m \in \mathbb{N}$  define
\ba
{\E}_{2p}^{m} & := & \{ R \in {\E}_{\infty} \; {\Big{|}} \;  (den(R))^{1/(m+1)} 
\text{ is even, symmetric of degree } 2p \}
\nn
\ea 
\no
and ${\E}_{2p}^{m,s} := {\E}_{2p}^{m} \cap {\E}_{2p}^{s}$, so a function
$R \in {\E}_{2p}^{m,s}$ can be written in the form 
\ba
R(z) & = &  \frac{P(z)}{Q^{m+1}(z) }, \nn
\ea
\no
where $P(z)$ is an even polynomial and $Q(z)$ is an even 
symmetric polynomial of degree  $2p$.

The method of partial fractions gives (in principle) the 
value of $I$ in terms of the roots of $Q$. Symbolic computations yield
either a closed-form answer, an expression in terms of the 
roots of an associated polynomial, or the integral returned unevaluated.

The algorithm described here allows only algebraic 
operations on elementary functions and changes of variables of the same type.
In particular, we exclude the solution of algebraic equations of degree higher 
than $2$.
We say that a rational function 
$R \in {\E}_{\infty}$ is {\em computable}  if its integral can be evaluated by 
our algorithm. 

The first step in the algorithm is to consider 
symmetric rational functions.  In Section \ref{sec-reduc}
we prove a reduction formula, i.e. a map 
${\mathfrak{F}}_{p} : {\E}_{2p}^{s} \to {\rat}_{p}$ that 
reduces the computability of the integral of the symmetric function 
$R$ to that of one
of degree $\tfrac{1}{2} deg(R)$. Here ${\rat}_{p}$ is the space of rational 
functions with denominator of degree $p$. These new
functions are not necessarily symmetric, and this is the main limitation 
of our algorithm. The details of 
${\mathfrak{F}}_{p}$ require the evaluation of some binomial sums which are 
presented in Appendix \ref{sec-binom}.  The classical Wallis' formula 
\ba
\ift \frac{dz}{(z^{2} + 1)^{m+1} } & = & \frac{\pi}{2^{2m+1}} \binom{2m}{m}
\nn
\ea
\no
shows that
every $R \in {\E}_{2}^{m}$ is computable. The reduction formula now 
implies that every $R \in {\E}_{4}^{m}$ is computable. This is 
described in Section
\ref{sec-quartic}.  The computability of $R \in {\E}_{4}$ is also a 
consequence of the 
classical theory of hypergeometric functions; the details are given
in \cite{hyp}. In Section \ref{sec-degree8} we compute the 
integral of every function in ${\E}_{8}^{m,s}$, where the reduction method 
expresses these integrals in terms of functions in ${\E}_{4}^{m}$.  The
algorithm does not, in general, provide a value for the integral of a 
nonsymmetric function of degree $8$. 

The final 
piece of the algorithm is described in Section \ref{sec-landen}: for 
$R \in {\E}_{2p}$, the 
symmetrization of its denominator produces a (symmetric) rational 
function in ${\E}_{4p}$ with the same integral as $R$. The
reduction formula now yields a new function in
${\E}_{2p}$ with the same integral as $R$. We thus obtain a map 
$ {\mathfrak{T}}_{2p}: {\E}_{2p} \to {\E}_{2p}$
such that 
\ba
\int_{0}^{\infty} R(z) \, dz & = & \int_{0}^{\infty} 
{{{\mathfrak{T}}_{2p}}(R(z))} \, dz.  \label{1invar}
\ea

\no
In particular, the class of computable rational functions of degree $2p$ is 
invariant under forward and backward iteration of ${\mathfrak{T}}_{2p}$. 
This map is the rational analog of the original Landen transformation 
for elliptic integrals described in \cite{borbor,mckmoll}.  The  map 
${\mathfrak{T}}_{2p}$ can also be interpreted as a map on the coefficients 
$\Phi_{2p} : {\mathfrak{O}}_{2p}^{+} \to 
{\mathfrak{O}}_{2p}^{+}$ where 
${\mathfrak{O}}_{2p}^{+} = 
{\mathbb{R}_{+}}^{p-1} \times {\mathbb{R}_{+}}^{p}$. 

The case of $\Phi_{6}$ is described in detail in  
Section \ref{sec-degree6} and
is given explicitly by 
\ba
a_{1} & \to & \frac{9 + 5a_{1} + 5a_{2} + a_{1}a_{2}}
{( a_{1} + a_{2} + 2)^{4/3} } \label{scheme6a}  \\
a_{2} & \to & \frac{a_{1} + a_{2}+6}{( a_{1} + a_{2} + 2)^{2/3} } \nn \\
b_{0} & \to & \frac{b_{0} + b_{1} + b_{2}}{(a_{1}+a_{2}+2)^{2/3} } \nn \\ 
b_{1} & \to & \frac{b_{0}(a_{2} + 2) + 2b_{1} + b_{2}(a_{1} + 3)}{ a_{1} + a_{2} + 2} \nn \\
b_{2} & \to & \frac{b_{0} + b_{2}}{( a_{1} + a_{2} + 2)^{1/3} }. \nn 
\ea

\no
Let ${\mathbf{x}}_{0}:= (a_{1},a_{2};b_{0},b_{1},b_{2})$. Then
$\Phi_{6}: {\mathfrak{O}^{+}_{6}} \to {\mathfrak{O}^{+}_{6}}$ is 
iterated to produce a sequence ${\mathbf{x}}_{n+1} := 
\Phi_{6}({\mathbf{x}}_{n})$
of points in ${\mathfrak{O}^{+}_{6}}$ that yield a sequence of
rational functions with constant integral. We have proved 
in \cite{landen} the existence of 
$L \in {\mathbb{R}^{+}}$, depending upon the initial 
point ${\mathbf{x}}_{0} = (a_{1},a_{2};b_{0},b_{1},b_{2}),$  such that  
${\mathbf{x}}_{n} \to 
(3,3;L,2L,L)$. Thus
\ba
\ift \frac{b_{0}z^{4} + b_{1}z^{2} + 
b_{2}}{z^{6} + a_{1}z^{4} + 
a_{2}z^{2} + 1} \; dz & = &  
L({\mathbf{x}}_{0}) \times \frac{\pi}{2}.  \label{int1}
\ea
\no
This establishes a numerical method to compute the integral in (\ref{int1}). 

Numerical studies on integrals of even degree $2p$ 
suggest the existence of a limiting value 
$L = L( {\mathbf{x}}_{0})$ such that the sequence 
${\mathbf{x}}_{n} := \Phi_{2p}( 
{\mathbf{x}}_{n} )$  satisfies 
\ba
{\mathbf{x}}_{n} & \to &  \left( 
\binom{p}{1}, \binom{p}{2}, \cdots, \binom{p}{p-1};
\binom{p-1}{0} \, L, \binom{p-1}{1} \, L, \cdots, 
\binom{p-1}{p-1} \, L \right).
\nn 
\ea
\no
The integral of the original rational function is thus
$\tfrac{\pi}{2} \times L({\mathbf{x}}_{0})$. The proof of convergence remains 
open for $p \geq 4$.  Examples are given in Section \ref{sec-examples}.

The most important
issues left unresolved in this paper are the convergence of the iteration 
of the map $\Phi_{2p}$ discussed above and the geometric interpretation of 
the Landen transformation ${\mathfrak{T}}_{2p}$. Finally, the question of
integration of odd functions has not been addressed at all.

\vskip 20pt

\no
{\bf Some history}. The problem of 
integration of rational functions $R(z)=P(z)/Q(z)$ was 
considered by 
J. Bernoulli in the $18^{th}$ century.  He completed the original 
attempt by Leibniz
of a general partial fraction decomposition of $R(z)$. 
The main difficulty associated with this procedure is to obtain a 
complete factorization of $Q(z)$ over $\R$. Once this is known
the partial fraction decomposition of $R(z)$ can be  computed. The 
fact is that the primitive of a rational function is always 
elementary: it consists of a new rational function (its {\em rational 
part}) and the logarithm of a second rational function (its {\em 
transcendental part}). In his classic monograph \cite{hardy} G. H. Hardy states:
{\em The solution of the problem (of definite integration) in the case of 
rational functions may therefore be said to be complete; for the difficulty 
with regard to the explicit solution of algebraical equations is one not of 
inadequate knowledge but of proved impossibility}. He goes on to add: 
{\em It appears from the preceding paragraphs that we can always find the
rational part of the integral, and can find the complete integral if we 
can find the roots of} $Q(z) = 0 $. 

In the middle of the last century Hermite \cite{herm} and Ostrogradsky
\cite{ostro} developed algorithms to compute the rational part of the primitive
of $R(z)$ {\em without} factoring $Q(z)$. More recently Horowitz \cite{horo}
rediscovered this method and discussed its complexity.  The problem of 
computing the transcendental part of the primitive was finally solved by 
Lazard and Rioboo \cite{lazrio}, Rothstein \cite{roth} and 
Trager \cite{trag}.  For detailed descriptions and proofs of these algorithms 
the reader is referred to \cite{bron} and \cite{geddes}.

\no
\section{The reduction formula} \label{sec-reduc}
\setcounter{equation}{0}

In this section we present a map ${\mathfrak{F}}_{p}: 
{\mathfrak{E}}_{2p}^{m,s} \to {\mathfrak{E}}_{p}^{m}$ that is the basis of the 
integration algorithm described in Section \ref{sec-landen}. The proof 
is elementary and the
binomial sums discussed in the Appendix are employed.  

Let $D_{p}(z)$ be the general symmetric polynomial of degree $4p$. We express
the integral of $z^{2n}/D_{p}^{m+1}$ as a linear combination of integrals of 
$z^{2j}/E_{p}^{m+1}$ where $E_{p}$ is a polynomial of degree $2p$ whose 
coefficients are determined by those of $D_{p}$.

\begin{Thm}
\label{main}
Let $m,n,p \in \mathbb{N}$. Define 
\ba
D_{p}(d_{1},d_{2}, \cdots, d_{p};z) & = & \sum_{k=0}^{p} d_{p+1-k} 
(z^{2k} + z^{4p-2k})  \label{dp}
\ea
\no
and 
\ba
 E_{p}(d_{1},d_{2}, \cdots, d_{p}; z) & = & 
 \left( \sum_{j=1}^{p+1} d_{j} \right) z^{2p} + \nn \\
 & + & \sum_{i=1}^{p} 2^{2i-1} z^{2(p-i)} \sum_{j=1}^{p-i+1} 
   \frac{j+i-1}{i} \binom{j+2i-2}{j-1} d_{j+i}, \nn \\ 
 & & \label{ep}
\ea
\no
for $d_{i} \in {\mathbb{R}}_{+}$. Then 
for $ \; 0 \leq n \leq (m+1)p-1,$ 
\ba
\int_{0}^{\infty} \frac{ z^{2n} \; dz }{ \left( D_{p}(d_{1}, \cdots,d_{p};z) 
\right)^{m+1}  }   =   \nn
\ea
\ba
2^{-m} 
\sum_{j=0}^{(m+1)p-n-1} 4^{j} \binom{(m+1)p-n-1+j}{2j} 
\int_{0}^{\infty} \frac{ 
z^{2((m+1)p-1-j)}} 
{ \left( \, E_{p}(d_{1},\cdots, d_{p};z) \, \right)^{m+1} } \; dz,  \nn \\
& & \label{newmain}
\ea
\no
and for $(m+1)p - 1 < n < 2p(m+1)-1$ we employ the symmetry rule 
\ba
N_{n,p} & = & N_{2p(m+1)-1-n,p}.  \label{symrule}
\ea
\end{Thm}

\no
\begin{proof}
First observe that (\ref{symrule}) follows from the change of variable
$z \to 1/z$. Now consider
\ba
N_{n,p}(d_{1}, \cdots, d_{p};m) & := & 
 \int_{0}^{\infty} \frac{ z^{2n} \; dz }
{ \left( \; \sum_{k=0}^{p} d_{p+1-k} (z^{4p-2k} + z^{2k} ) \; \right)^{m+1}}
\label{nnp}
\ea
\no
for $0 \leq n \leq (m+1)p-1$. The 
substitution $z = \tan \, \theta$ yields 
\ba
N_{n,p} & = & 
   \int_{0}^{\pi/2}  \frac{ (1 - C^{2})^{n} C^{4(m+1)p-2n-2} \; d \theta }
{ \left( \; \sum_{k=0}^{p} d_{p+1-k} \left\{ (1-C^{2})^{2p-k} C^{2k} + 
(1-C^{2})^{k} C^{4p-2k} 
\right\} \; \right)^{m+1} },  \nn 
\ea
\no
where $C = \cos \theta$. Letting 
$\psi = 2 \, \theta$ and
$D = \cos \psi = 2C^{2} -1$ then gives 
\ba
N_{n,p} & = & \int_{0}^{\pi} 
\frac{ (1 - D^{2})^{n} (1+D)^{2(m+1)p-2n-1} \; d \psi }
{ \left( \; \sum_{k=0}^{p} d_{p+1-k} (1-D^{2})^{k} \left\{(1-D)^{2p-2k} + 
(1+D)^{2p-2k}  \right\}  \; \right)^{m+1} }. \nn
\ea
\no
Now observe that the integrals of the odd powers of cosine vanish when we 
expand $(1+D)^{2(m+1)p-2n-1}$, producing
\ba
N_{n,p} & = & 2^{-m} \int_{0}^{\pi/2}
\frac{(1-D^{2})^{n} \sum_{j=0}^{(m+1)p-n-1} \binom{2(m+1)p-2n-1}{2j} 
D^{2j} \, d \theta}
{ \left\{ \sum_{k=0}^{p} d_{p+1-k} (1-D^{2})^{k} \sum_{j=0}^{p-k} 
\binom{2p-2k}{2j} D^{2j} \right\}^{m+1} }.  \nn 
\ea
\no
A second double angle substitution $\varphi = 2 \psi$  gives 
\ba
N_{n,p} & = & 2^{-m} \int_{0}^{\pi}
\frac{(1-E)^{n} \sum_{j=0}^{(m+1)p-n-1} \binom{2(m+1)p-2n-1}{2j} 
2^{p(m+1)-n-j-1} (1+E)^{j} \, d \varphi}
{ \left\{ \sum_{k=0}^{p} d_{p+1-k} (1-E)^{k} \sum_{j=0}^{p-k} 
\binom{2p-2k}{2j} 2^{p-k-j} (1+E)^{j} \right\}^{m+1} },  \nn
\ea
\no
where $E  = \cos \varphi = 2D^{2} - 1$. The 
change of variable $z = \tan(\varphi/2)$ then yields 
\ba 
N_{n,p} & = & 2^{-m} \int_{0}^{\infty} 
\frac{ z^{2n}  \sum_{j=0}^{(m+1)p-n-1} \binom{2(m+1)p-2n-1}{2j} 
(1+z^{2})^{(m+1)p-n-j-1}
\, dz } 
{ \left\{ \sum_{k=0}^{p} d_{p+1-k} z^{2k} \left( \sum_{j=0}^{p-k} 
\binom{2p-2k}{2j} (1+z^{2})^{p-k-j} \right) \right\}^{m+1} }. \nn  \\
 & & \label{mess10}
\ea
\no
Finally, we modify (\ref{mess10}) using 
Lemma \ref{lem32} and Lemma \ref{lem34} with 
$N = (m+1)p-n-1$ to produce  (\ref{newmain}).
\end{proof}

Note that the previous theorem associates to each rational function 
of symmetric denominator
\ba
R_{1}(z) & = & \frac{b_{s}z^{2s} + b_{s-1}z^{2(s-1)} + \cdots + b_{0} }
{ \left( \, z^{4p} + d_{p} z^{4p-2} + \cdots + 2d_{1}z^{2p} + \cdots 
+ 1 \, \right)^{m+1} }  \nn
\ea
\no
a new rational function
\ba
R_{2}(z) & = & 2^{-m} \sum_{n=0}^{(m+1)p-1}  
b_{n} \sum_{j=0}^{(m+1)p-n-1} 4^{j} \binom{(m+1)p-n-1+j}{2j} 
\frac{z^{2((m+1)p-1-j)}} 
{ \left( E_{p}(d_{1},\cdots, d_{p};z) \right)^{m+1} } \; \nn
\ea
\no
such that 
\ba
\ift R_{1}(z) \, dz & = & \ift R_{2}(z) \, dz. \nn
\ea

\vskip 5pt

\no
\section{The quartic case} \label{sec-quartic}
\setcounter{equation}{0}

In this section we describe the computability of rational functions 
$R \in {\E}_{4}^{m}$. These are functions of the form 
\ba
R(z) & = &  \frac{P(z)}{(z^{4} + 2az^{2} + 1)^{m+1} } \nn
\ea
\no
where $P(z)$ is an even polynomial of degree $4m+2$. Observe that
the normalization $a_{0}=a_{2} = 1$ makes the denominator of $R$ 
automatically symmetric.  It suffices to evaluate
\ba
N_{n,4}(d_{1};m) & := &  \int_{0}^{\infty} \frac{z^{2n} \; dz}{(z^{4} + 2d_{1}z^{2} 
+1)^{m+1} } \label{n04}
\ea
\no
where $0 \leq n \leq 2m+1$ is required for convergence. From (\ref{symrule}) 
we have
$N_{n,4}(d_{1};m)  =  N_{2m-1-n,4}(d_{1};m)$, so we
may assume $0 \leq n \leq m$. 
We now employ Theorem \ref{main} to obtain a closed form expression 
for $N_{n,4}(d_{1};m)$.  

\begin{Thm}
\label{mainthm}
Let $m \in \mathbb{N}$ and assume $0 \leq n \leq m$. Then 
\ba
N_{n,4}(d_{1};m) := \int_{0}^{\infty} \frac{z^{2n} \; dz}
{ \left( z^{4} + 2d_{1}z^{2} + 1 \right)^{m+1} } \; \; =  \label{quartic1}
\ea
\ba
\frac{\pi}{ 2^{3m+3/2} (1+a)^{m+1/2} } \times 
\sum_{j=0}^{m-n} 2^{j}(1+d_{1})^{j} \times
\binom{2m-2j-1}{m-j}  \binom{m-n+j}{2j} \binom{2j}{j}
\binom{m}{j}^{-1} . \nn
\ea
\no
For $m+1 \leq n \leq 2m+1$ we have
\ba
\int_{0}^{\infty} \frac{z^{2n} \; dz}
{ \left( z^{4} + 2d_{1}z^{2} + 1 \right)^{m+1} } \; \; =  \label{quartic2}
\ea
\ba
\frac{\pi}{ 2^{3m+3/2} (1+d_{1})^{m+1/2} } \times 
\sum_{j=0}^{n-m-1} 2^{j}(1+d_{1})^{j} \times
\binom{2m-2j-1}{m-j}  \binom{m-n+j}{2j} \binom{2j}{j}
\binom{m}{j}^{-1} . \nn
\ea
\end{Thm}

\begin{proof}
We apply the result of the Theorem \ref{main} with 
$D_{1}(d_{1};z)  =  z^{4} + 2d_{1}z^{2} + 1 $ and
$E_{1}(d_{1};z) = (1 + d_{1})z^{2} + 2$, so that 
\ba
\int_{0}^{\infty} \frac{z^{2n} \, dz }{(z^{4} + 2d_{1}z^{2} + 1)^{m+1} } & = &
2^{-m} \sum_{j=0}^{m-n} 4^{j} \binom{m-n+j}{2j} 
\int_{0}^{\infty} \frac{z^{2(m-j)} \, dz}{ ( (1+d_{1})z^{2} + 2 )^{m+1} }. \nn
\ea
\no
The change of variable $u = (1+d_{1})z^{2}/2$ then yields
\ba
\int_{0}^{\infty} \frac{z^{2(m-j)} \, dz}{ ( \; (1+d_{1})z^{2} + 2 \; )^{m+1} } 
& = &  2^{-(j+3/2)} (1+d_{1})^{-m+j-1/2} 
\int_{0}^{\infty} \frac{u^{m-j-1/2} \; du }{(1+u)^{m+1} }
\nn \\
 &  = & \pi \binom{2m-2j}{m-j} \binom{2j}{j} \binom{m}{j}^{-1} 2^{-(2m+j+3/2)}
(1+d_{1})^{-(m-j+1/2)}, \nn
\ea
\no
where we have used
\ba
\ift \frac{u^{r- 1/2}}{(1+u)^{s}} \, du & = & 
\frac{\pi}{2^{2(s-1)}} \binom{2r}{r} \binom{2(s-r-1)}{s-r-1} 
\binom{s-1}{r}^{-1}. \nn
\ea
\end{proof}

\no
The algorithm also requires a scaled version of $N_{0,4}(d_{1};m)$. 

\no
\begin{Cor}
\label{cor62}
Let $b>0, \, c>0,\; a > - \sqrt{bc}, \, m \in \mathbb{N}$, 
and $0 \leq n \leq m$. Define
\ba
N_{n,4}(a,b,c;m) & := &  \ift \frac{z^{2n} \; dz}{ \left( bz^{4} + 2az^{2} + c 
\right)^{m+1}}. \nn
\ea
\no
Then for $0 \leq n \leq m$, 
\ba
N_{n,4}(a,b,c;m) & = & 
 \pi \left( c (c/b)^{m-n} \left\{ 8(a + \sqrt{bc}) \right\}^{2m+1} 
\right)^{-1/2} \times  \nn \\
 & & \times
\sum_{k=0}^{m-n} 2^{k} \binom{2m-2k}{m-k} \binom{m-n+k}{2k} 
\binom{2k}{k} \binom{m}{k}^{-1} \left( \frac{a}{\sqrt{bc}} +1 
\right)^{k}, \nn \\
 &  & \label{sum1}
\ea
and for $m+1 \leq n \leq 2m+1$, 
\ba
N_{n,4}(a,b,c;m) & = & 
 \pi \left( c (c/b)^{m-n} \left\{ 8(a + \sqrt{bc}) \right\}^{2m+1} 
\right)^{-1/2} \times  \nn \\
 & & \times
 \sum_{k=0}^{n-m-1} 2^{k} 
\binom{2m-2k}{m-k} \binom{m-n+k}{2k} 
\binom{2k}{k} \binom{m}{k}^{-1} \left( \frac{a}{\sqrt{bc}} + 1 \right)^{k}.
\nn \\
 &  & \label{sum2}
\ea
\end{Cor}

\begin{proof}
Let $0 \leq n \leq m$. The substitution $u = z (b/c)^{1/4}$ yields
\ba
N_{n,4}(a,b,c;m) & = & 
\frac{1}{c^{m-n/2 + 3/4} b^{n/2 + 1/4} } N_{n,4} 
\left( \frac{a}{\sqrt{bc}} ; m \right), \label{scaled}
\ea
\no
so (\ref{sum1}) then follows from Theorem \ref{mainthm}.
From (\ref{symrule}) we have
$N_{n,4}(a,b,c;m)  =  
N_{2m+1-n,4}(c,b,a;m)$ for $m+1 \leq n \leq 2m+1$, giving (\ref{sum2}).
\end{proof}

\section{The symmetric case of degree $8$} \label{sec-degree8}
\setcounter{equation}{0}

In this section we prove the computability of the set ${\E}_{8}^{m,s}$ of 
symmetric rational functions with denominator of degree $8$ and establish an 
explicit formula for the integral 
\ba
N_{n,8}(a_{1},a_{2};m) & = & \int_{0}^{\infty} 
\frac{z^{2n} \; dz} { (z^{8} + a_{2} z^{6} + 2a_{1}z^{4} + a_{2}z^{2} 
+ 1)^{m+1} } \nn
\ea
\no
where $0 \leq n \leq 4m+3$ is required for convergence. Observe that 
(\ref{symrule}) reduces the discussion to the case  $0 \leq n \leq 2m+1$. 
The expression (\ref{ep}), with $p=2$,  produces
$E_{2}(a_{1},a_{2};z) = (1 + a_{1} + a_{2}) z^{4} + 2(a_{2} + 4)z^{2} + 8$. 

\begin{Thm}
\label{thmsym8}
Every function in ${\E}_{8}^{m,s}$ is computable. More specifically, define
$c_{1} := a_{2} + 4, \; c_{2} := 1 + a_{1} + a_{2}$, and
\ba
t_{k,j}(m,n;a_{1},a_{2}) & := &  \pi  2^{-(3m+2+k+j)/2} 
c_{2}^{(m-k-j)/2} (c_{1} + \sqrt{8c_{2}})^{j-m-1/2} \times  \nn \\
 & \times & \binom{4m-n-k+2}{k-n} \binom{2m-2j}{m-j} \binom{m-k+j}{2j}
\binom{2j}{j} \binom{m}{j}^{-1}. \nn
\ea
\no
Then for $m+1 \leq n \leq 2m+1, \; 1+a_{1}+a_{2} > 0$ and 
$a_{2}+4 > -8\sqrt{8(1+a_{1}+a_{2})}$, 
\ba
\int_{0}^{\infty} \frac{z^{2n} \; dz}{(z^{8} + a_{2}z^{6} + 2a_{1}z^{4} +
a_{2} z^{2} + 1)^{m+1} } & = & 
\sum_{k=n}^{2m+1} \sum_{j=0}^{k-m-1} t_{k,j}(m,n;a_{1},a_{2}), \nn
\ea
\no
and for $0 \leq n \leq m$, 
\ba
\int_{0}^{\infty} \frac{z^{2n} \; dz}{(z^{8} + a_{2}z^{6} + 2a_{1}z^{4} +
a_{2} z^{2} + 1)^{m+1} }  = \nn \\
\sum_{k=n}^{m} \sum_{j=0}^{m-k} t_{k,j}(m,n;a_{1},a_{2}) +
\sum_{k=m+1}^{2m+1} \sum_{j=0}^{k-m-1} t_{k,j}(m,n;a_{1},a_{2}). \nn
\ea
\end{Thm}
\begin{proof}
The reduction formula yields 
\ba
\int_{0}^{\infty} \frac{z^{2n} \; dz}{(z^{8} + a_{2}z^{6} + 2a_{1}z^{4} +
a_{2} z^{2} + 1)^{m+1} }  = \nn \\
2^{3m+2} \sum_{k=n}^{2m+1} 2^{-2k} \binom{4m-n-k+2}{k-n} 
\int_{0}^{\infty} \frac{ z^{2k} \; dz}{ (c_{2} z^{4} + 2c_{1} z^{2} + 8)^{m+1} 
}. \label{label1}
\ea
\no
We then use Corollary \ref{cor62} to evaluate (\ref{label1}). 
\end{proof}

\section{A sequence of Landen transformations} \label{sec-landen}
\setcounter{equation}{0}

The transformation theory of elliptic integrals was initiated 
by Landen in 1771. He proved the invariance of the function
\ba
G(a,b) & := & \int_{0}^{\pi/2} \frac{d \, \theta} 
{ \sqrt{ a^{2} \cos^{2} \theta + b^{2} \sin^{2}  \theta } }   \label{gab}
\ea
\no
under the transformation 
\ba
a_{1} = (a+b)/2 &  \qquad &  b_{1} = \sqrt{ab}, \label{transgl}
\ea
\no
i.e. that
\ba
G(a_{1},b_{1}) & = & G(a,b). \label{invargl}
\ea
\no
Gauss \cite{gauss1} rediscovered this invariance 
while numerically calculating the length of 
a lemniscate.  
An elegant proof of (\ref{invargl}) is given by 
Newman in \cite{newman}. Here, the substitution $x = b \tan \theta$ 
converts $2 G(a,b)$ into the integral of $\left[ (a^{2} + x^{2})
(b^{2} + x^{2}) \right]^{-1/2}$ over $\mathbb{R}$; the change of variable 
$t = (x-ab/x)/2$ then completes the proof. 

The Gauss-Landen transformation can be iterated to produce a 
double sequence $(a_{n},b_{n})$ such that $0 \leq a_{n} - b_{n} < 2^{-n}$.
It follows that $a_{n}$ and $b_{n}$ converge to a common limit, the so-called 
{\em arithmetic-geometric mean} of $a$ and $b$, denoted by $AGM(a,b)$. Passing 
to the limit in $G(a,b) = G(a_{n},b_{n})$ produces 
\ba
\frac{\pi}{2 AGM(a,b) } & = &
\int_{0}^{\pi/2} \frac{d \, \theta} 
{ \sqrt{ a^{2} \cos^{2} \theta + b^{2} \sin^{2}  \theta }}.
\ea
\no
The reader is referred to \cite{borbor} and \cite{mckmoll} for details.

The goal of this section is to produce a map 
${{\mathfrak{T}}}_{2p}: {{\mathfrak{E}}_{2p}} \to 
{{\mathfrak{E}}_{2p}}$
that preserves the integral, i.e.
\ba
\int_{0}^{\infty} R(z) \, dz & = & \int_{0}^{\infty} 
{{{\mathfrak{T}}_{2p}}(R(z))} \, dz.  \label{invar}
\ea
\no
This map is the rational analog of the original Landen transformation 
(\ref{transgl}).

\begin{Thm}
\label{thmbig}
Let $R(z) = P(z)/Q(z)$ with
\ba
P(z) = \sum_{j=0}^{p-1} b_{j} z^{2(p-1-j)} \text{   and   } \; \;
Q(z) = \sum_{j=0}^{p} a_{j} z^{2(p-j)}. \label{polyb}
\ea
\no
Define $a_{j}=0$ for $j>p$, $b_{j}=0$ for $j>p-1$, 
\ba
d_{p+1-j} & = & 
\sum_{k=0}^{j} a_{p-k} a_{j-k} 
\label{dofj}
\ea
\no
for $0 \leq k \leq p-1$,
\ba
d_{1} & = & 
\frac{1}{2} \sum_{k=0}^{p} a_{p-k}^{2},
\label{dofj1}
\ea
\ba
c_{j} & = & \sum_{k=0}^{2p-1} a_{j}b_{p-1-j+k} 
\label{cofj}
\ea
\no
for $0 \leq j \leq  2p-1$, and also
\ba
\alpha_{p}(i) & = & \begin{cases}
     2^{2i-1} \sum_{k=1}^{p+1-i} \frac{k+i-1}{i} \binom{k+2i-2}{k - 1}
d_{k + i} \text{    if } 1 \leq i \leq p    \\
     1 + \sum_{k=1}^{p} d_{k} \text{ if } i = 0. 
          \end{cases} 
\label{alpha}
\ea
\no
Let
\ba
a_{i}^{+} & = & \frac{\alpha_{p}(i)}{2^{2i} Q(1)^{2(1-i/p)}} 
\label{acoeff}
\ea
\no
for $1 \leq i \leq p-1$, and
\ba
b_{i}^{+} & = & Q(1)^{2i/p+1/p-2}  \times 
\left[ \sum_{k=0}^{p-1-i} (c_{k} + c_{2p-1-k}) \binom{p-1-k+i}{2i} 
\right] \label{bcoeff}  
\ea
\no
for $0 \leq i \leq p-1$. 
\no
Finally, define the polynomials 
\ba
P^{+}(z)  = {\sum_{k=0}^{p-1} b_{i}^{+} z^{2(p-1-i)}} & \text{ and } &
Q^{+}(z) = {\sum_{k=0}^{p} a_{i}^{+} z^{2(p-i)}}. \label{newb}
\ea
\no
Then ${\mathfrak{T}}_{2p}(R(z)) := P^{+}(z)/Q^{+}(z)$ satisfies 
{\rm{(\ref{invar})}}, i.e.
\ba
\int_{0}^{\infty} \frac{P(z)}{Q(z)} \; dz & = & 
\int_{0}^{\infty} \frac{P^{+}(z)}{Q^{+}(z)} \; dz. 
\ea
\end{Thm}

\begin{proof}
The first step is to convert the polynomial $Q(z)$ to its symmetric form: 
\ba
I := \int_{0}^{\infty} \frac{P(z)}{Q(z)} dz  = 
\int_{0}^{\infty} \frac{C(z)}{D(z)} dz \nn
\ea
\no
with 
\ba
C(z) & = & P(z) \times z^{2p} Q(1/z) := \sum_{k=0}^{2p-1} c_{k}z^{2k} \nn 
\\
D(z) & = & Q(z) \times z^{2p} Q(1/z) := \sum_{k=0}^{p} d_{p+1-k}
(z^{2k} + z^{2(2p-k)} ). \nn
\ea
\no
Then
\ba
I  & =  & \sum_{k=0}^{2p-1} c_{k} \int_{0}^{\infty} \frac{z^{2k} \; dz}
{Q^{s}(z)}. \nn
\ea
\no
Now employ the reduction formula in Section \ref{sec-reduc} to evaluate 
\ba
L_{k} & := & \int_{0}^{\infty} \frac{z^{2k} \; dz }{Q^{s}(z)}. \nn
\ea
\no
Observe that one needs to evaluate $L_{k}$ only for $0 \leq k \leq p-1$. 
Indeed, the usual symmetry rule yields 
$L_{k} = L_{2p-1-k}$.  The reduction formula now gives
\ba
L_{k} & = & 
\sum_{j=0}^{p-1-k} 2^{2j} \binom{p-1-k+j}{2j}  
\int_{0}^{\infty} \frac{z^{2(p-1-j)} \; dz}{ \sum_{i=0}^{p}
\alpha_{p}(i) z^{2(p-i)} } \nn \\
& = & \frac{1}{\alpha_{p}(p)} \sum_{j=1}^{p-k} 
2^{2(j-1)} \binom{p-k+j-2}{2j-2} \lambda^{2p-2j+1}  
\int_{0}^{\infty} \frac{z^{2(p-j)} \; dz}{\sum_{i=0}^{p} b_{i}^{+} z^{2(p-i)} }
\nn
\ea
\no
with $\alpha_{p}(i)$ as in (\ref{alpha}) and
$\lambda = \left[ \alpha_{p}(p)/\alpha_{p}(0) \right]^{1/2p}$. 
\end{proof}

\no
{\bf Note}. The extension of this transformation to the case  of 
\ba
\ift \frac{P(z)}{Q^{m+1}(z)} \, dz \nn
\ea
\no
requires explicit formulae for the coefficients of $P(z) \times 
\left( \, z^{2p} Q(1/z) \right)^{m+1}$ and 
$ Q^{m+1}(z) \times \left( \, z^{2p} Q(1/z) \right)^{m+1}$.
\vskip 10pt

\no
{\bf An algorithm for integration}.  Let ${\mathbf{x}} = 
( {\mathbf{a}}, {\mathbf{b}})$ with ${\mathbf{a}} = (a_{1}, \cdots, a_{p-1})$,
${\mathbf{b}} = (b_{0}, \cdots, b_{p-1})$, and let ${\mathfrak{O}}_{2p}^{+} = 
{\mathbb{R}_{+}}^{p-1} \times {\mathbb{R}_{+}}^{p}$. We then have a 
map 
\ba
\Phi_{2p}: {\mathfrak{O}}_{2p}^{+} & \to & {\mathfrak{O}}_{2p}^{+}   \nn \\
     {\mathbf{x}}:= ( {\mathbf{a}}, {\mathbf{b}} ) & \to & 
     {\mathbf{x}}^{+}:= ( {\mathbf{a}}^{+}, {\mathbf{b}}^{+} )  \nn
\ea
\no
where $a_{i}^{+}$ and $b_{i}^{+}$ are given in (\ref{acoeff},
\ref{bcoeff}).  
Iteration of this map, starting at ${\mathbf{x}}_{0}$, produces a sequence  
${\mathbf{x}}_{n+1} := \Phi_{2p}({\mathbf{x}}_{n})$ of points in 
${\mathfrak{O}}_{2p}^{+}$. The rational functions formed with these parameters
have integrals that remain constant along this orbit. Numerical studies
suggest the existence of a number $L = L({\mathbf{x}}_{0}) 
\in {\mathbb{R}}_{+}$ such that
\ba
{\mathbf{x}}_{n} & \to &  \left( 
\binom{p}{1}, \binom{p}{2}, \cdots, \binom{p}{p-1};
\binom{p-1}{0} \, L, \binom{p-1}{1} \, L, \cdots, 
\binom{p-1}{p-1} \, L \right).
\nn 
\ea
\no
Thus the integral of the original rational function is
$\tfrac{\pi}{2} \times L$. \\

\section{The sixth degree case} \label{sec-degree6}
\setcounter{equation}{0}

We discuss the map 
${\mathfrak{T}}_{2p}: {{\mathfrak{E}}_{2p}} \to 
{{\mathfrak{E}}_{2p}}$
for the case $p=3$. The effect of ${\mathfrak{T}}_{6}$ on the 
coefficients 
${\mathfrak{P}}_{6} =  \{ b_{0},b_{1},b_{2},a_{1},a_{2} \}$
is denoted by $\Phi_{6} : {{\mathfrak{O}}}_{6}^{+} \to 
{{\mathfrak{O}}}_{6}^{+}$ and is given explicitly by 

\ba
a_{1} & \to & \frac{9 + 5a_{1} + 5a_{2} + a_{1}a_{2}}
{( a_{1} + a_{2} + 2)^{4/3} } \label{scheme6} \\
a_{2} & \to & \frac{a_{1} + a_{2}+6}{( a_{1} + a_{2} + 2)^{2/3} } \nn \\
b_{0} & \to & \frac{b_{0} + b_{1} + b_{2}}{(a_{1}+a_{2}+2)^{2/3} }  \nn \\
b_{1} & \to & \frac{b_{0}(a_{2} + 2) + 2b_{1} + b_{2}(a_{1} + 3)}{ a_{1} + a_{2} + 2} \nn \\
b_{2} & \to & \frac{b_{0} + b_{2}}{( a_{1} + a_{2} + 2)^{1/3} } \nn
\ea
\no
using Theorem \ref{thmbig}. The 
map $\Phi_{6}$  preserves the integral 
\ba
U_{6}(a_{1},a_{2},b_{0};b_{1},b_{2}) & := & 
\ift \frac{b_{0}z^{4} + b_{1}z^{2} + b_{2}}{z^{6} + a_{1}z^{4} + 
a_{2}x^{2} + 1} \; dz
\ea
\no
and the convergence of its iterations has been proved in \cite{landen}, the
main result of which is the following theorem. 

\begin{Thm}
Let ${\mathbf{x}}_{0}:= (a_{1}^{0},a_{2}^{0};b_{0}^{0},b_{1}^{0},b_{2}^{0})
\in {\mathbb{R}_{+}}^{5}$. Define
${\mathbf{x}}_{n+1} := \Phi_{6}({\mathbf{x}}_{n})$.  Then 
$U_{6}$ is invariant under $\Phi_{6}$.  Moreover, the 
sequence $\{ (a_{1}^{n},a_{2}^{n}) \}$ converges to $(3,3)$ and
$\{ (b_{0}^{n}, b_{1}^{n},b_{2}^{n}) \}$ converges 
to $(L, 2L,L)$, where the limit 
$L$ is a function of the initial data ${\mathbf{x}}_{0}$. 
Therefore
\ba
\ift \frac{b_{0}z^{4} + b_{1}z^{2} + 
b_{2}}{z^{6} + a_{1}z^{4} + 
a_{2}z^{2} + 1} \; dz & = &  
L({\mathbf{x}}_{0}) \times \frac{\pi}{2}.  \nn
\ea
\end{Thm} 

\no
This iteration is similar to Landen's transformation for elliptic integrals
that has been employed in \cite{borbor} in the efficient calculation of $\pi$. 
Numerical data indicate
that the convergence of ${{\mathbf{x}}}_{n}$ is 
quadratic. The proof of convergence is based  on the fact that $\Phi_{6}$
cuts the distance from $(a_{1},a_{2})$ to $(3,3)$ by at least half.  \\

\no
{\bf A sequence of algebraic curves}. The complete characterization of 
parameters $(a_{1},a_{2})$ in the first quadrant  that yield computable 
rational functions 
\ba
R(z) & := & 
\frac{b_{0}z^{4} + b_{1}z^{2} + 
b_{2}}{z^{6} + a_{1}z^{4} + 
a_{2}z^{2} + 1}. \nn
\ea
\no
of degree $6$ remains open. The  polynomial 
$z^{6} + a_{1}z^{4} + a_{2}z^{2} + 1$ factors when 
$a_{1} = a_{2}$ so the diagonal
$\Delta  :=  \{ (a_{1},a_{2}) \in  {\R}_{+} \times {\R}_{+}: a_{1} = a_{2} \} $
produces computable functions. 
In view of the invariance of the class of computable functions 
under iterations by $\Phi_{6}$, the  curves 
${\mathbb{X}}_{n} := \Phi_{6}^{(-n)}(\Delta)$, with $n \in {\mathbb{Z}}$, are 
also computable.

The curve  ${\mathbb{X}}_{1}$ has equation 
\ba
(9 + 5a_{1} + 5a_{2} + a_{1}a_{2})^{3} & = & (a_{1} + a_{2} + 2)^{2} 
(a_{1} + a_{2} + 6)^{3}  \nn
\ea
\no
and consists of two branches meeting at the cusp $(3,3)$. In terms 
of the coordinates $x = a_{1} - 3$ and $y = a_{2} - 3$ the leading order 
term is $T_{1}(x,y) = 1728(x-y)^2$. This curve 
is rational and can be parametrized by 
\ba
a_{1}(t) & = & t^{-2}( t^{5} -t^{4} + 2t^{3} -t^{2} + t + 1) \label{param1} \\
a_{2}(t) & = & t^{-3}( t^{5} +t^{4} - t^{3} + 2t^{2} - t + 1). \nn
\ea

The rationality of ${\mathbb{X}}_{n}$ for $n \neq 1$ and its significance 
for the integration algorithm remains open. The  complexity 
of these curves increases with $n$. For example, the curve 
${\mathbb{X}}_{2} := \Phi_{3}^{(-2)}(\Delta)$ is of total degree $90$ in 
$x = a_{1} -3$ and $y=a_{2}-3$ with leading term
\ba
T_{2}(x,y) & := & 2^{121} 3^{35} (x-y)^{18} \left[ -163(x^{4} + y^{4}) +
668xy(x^{2} + y^{2}) -1074x^{2}y^{2} \right]. \nn
\ea
\no
The diagonal $\Delta$ can be replaced by a $2$-parameter family of 
computable curves ${\mathbb{X}}(c,d)$ that are produced from the factorization 
of the sextic with $a_{1} = c+d$ and $a_{2} = cd + 1/d$.  All the images 
$\Phi_{6}^{(-n)}{\mathbb{X}}(c,d)$ with $n \in {\mathbb{Z}}$ are computable
curves.  The question of whether these are {\em all }
the computable parameters remains open.

\section{Examples} \label{sec-examples}
\setcounter{equation}{0}

In this section we present a variety of closed-form evaluations of integrals of 
rational functions. \\

\no
{\bf Example 1}. The integral 
\ba
\ift \frac{z^2}{(z^4 + 4z^2 + 1)^9} \; dz & = & \frac{23698523 \, \pi}
{12230590464 \sqrt{6} } \nn
\ea
\no
is computed by Mathematica 3.0 using (\ref{quartic1}) in $.01$ seconds. The 
direct calculation took $12.27$ seconds and $6.4$ extra seconds to simplify 
the answer. \\

\no
{\bf Example 2}. The integral of any even rational function with denominator a 
power of an even quartic polynomial can be computed directly by using 
Corollary \ref{cor62}. For example:
\ba
\ift \frac{z^{6} \, dz}{ (2z^{4} + 2z^{2} + 3)^{11}  } & = & \frac{11 \pi  
(14229567 + 4937288 \, \sqrt{6})}{440301256704 \, ( 1 + \sqrt{6} \, )^{21/2} }.
\nn
\ea
\vskip 10pt

\no
{\bf Example 3}. The case $n=0$ in (\ref{quartic1}) deserves special attention:
\ba
N_{0,4}(a;m) & = &  \frac{\pi}{2^{m+3/2} (a+1)^{m+1/2}} P_{m}(a) \label{no4}
\ea
\no
where 
\ba
P_{m}(a) & = & 2^{-2m} \sum_{k=0}^{m} 2^{k} \binom{2m-2k}{m-k} \binom{m+k}{m} 
(a+1)^{k}. \label{polyam}
\ea
\no
The polynomial $P_{m}(a)$ has been studied in \cite{unim} and \cite{hyp}. \\

\no
{\bf Example 4}. The case $n = m$ in (\ref{quartic1}) yields 
\ba
N_{m,4}(a;m) =  
\int_{0}^{\infty} \frac{z^{2m} \; dz}{(z^{4} + 2az^{2} + 
1)^{m+1}} 
  =  \frac{\pi}{2^{3m+3/2} (1+a)^{m+1/2}} \times \binom{2m}{m}. \nn
\ea
\no
The change of variable $z \to \sqrt{z}$ converts this integral to 
\ba
N_{m,4}(a;m) & = & \frac{1}{2} \int_{0}^{\infty} 
\frac{z^{m-1/2} \; dz}{(z^{2} + 2az + 1)^{m+1} }, \nn
\ea
\no
which is  \cite{gr} 3.257.9.

\vskip 10pt

\no
{\bf Example 5}. A symmetric function of degree $6$. 
The integral
\ba
I  =   \ift \frac{x^{8}}{(x^{6} + 4x^{4} + 4x^{2} + 1)^{5}} \, dx  = 
\ift  \frac{x^{8}}{\left[ (x^{2} + 1)(x^{4} + 3x^{2} + 1 ) 
\right]^{5} } \, dx \nn
\ea
\no
can be computed by decomposing the integrand into partial fractions as
\ba
- \frac{1}{(x^{2} + 1)^{5}} - \frac{1}{(x^{2} + 1)^{4}} 
- \frac{6}{(x^{2} + 1)^{3}} - \frac{11}{(x^{2} + 1)^{2}} 
 -   \frac{31}{(x^{2} + 1)} \nn \\
+ \frac{1}{(x^{4} + 3x^{2} +1)^{5}} + 
\frac{2x^{2}}{(x^{4} + 3x^{2} + 1)^{5}} 
 -   \frac{4}{(x^{4} + 3x^{2} +1)^{4}} 
-\frac{3x^{2}}{(x^{4} + 3x^{2} + 1)^{4}}  + \nn \\
\frac{12}{(x^{4} + 3x^{2} +1)^{3}} +
   \frac{6x^{2}}{(x^{4} + 3x^{2} + 1)^{3}} 
 - \frac{32}{(x^{4} + 3x^{2} +1)^{2}}
 - \frac{14x^{2}}{(x^{4} + 3x^{2} + 1)^{2}}  \nn \\
  +   \frac{73}{(x^{4} + 3x^{2} +1)} + 
  \frac{31x^{2}}{(x^{4} + 3x^{2} + 1)}.  \nn
\ea
\no
Each of these terms is now computable yielding
\ba
I & = & \frac{1407326 \sqrt{5} - 3146875}{160000} \times \pi. \nn
\ea

\vskip 10pt

\no
{\bf Example 6}. Non-symmetric functions of degree $6$. In this case we 
can use the scheme (\ref{scheme6}) to produce 
numerical approximations to the integral. For example, the evaluation of
\ba
\ift \frac{45z^{4} + 25000z^{2} + 1230}{z^{6} + z^{4} + 3000z^{2} + 1} \; dz
 \nn
\ea
\no
is shown below:

\begin{center}
\begin{tabular}[h]{||c|c|c|c|c|c||}
\hline
$n$ & $a_{1}^{n}$ & $a_{2}^{n}$ & $b_{0}^{n}$ & $b_{1}^{n}$ & $b_{2}^{n}$  
\\ \hline 
0 & 1 & 3000 & 45 & 25000 & 1230  \\
1 &.415786 & 14.4465 & 126.233 & 63.2884 & 88.3741 \\
2 & 2.06562 & 3.17262 & 42.2607 & 156.015 & 83.6896 \\
3 & 2.98142 & 3.00338 & 75.3541 & 137.717 & 65.1111 \\
4 & 2.99999 & 3. & 69.6338 & 139.925 & 70.2771 \\
5 & 3. & 3. & 69.9589 & 139.914 & 69.9555 \\
6 & 3. & 3. & 69.9572 & 139.914 & 69.9572 \\
7 & 3. & 3. & 69.9572 & 139.914 & 69.9572 \\
\hline
\end{tabular}
\end{center}

\vskip 10pt

\no
Thus $L \sim 69.9572$ and 
\ba
\ift \frac{45x^{4} + 25000x^{2} + 1230}{x^{6} + x^{4} + 3000x^{2} + 1}
 \; dx & \sim &
69.9572 \times \frac{\pi}{2} \cong 109.889. \nn
\ea

\vskip 5pt 

\no
{\bf Example 7}. Symmetric functions of degree $8$. These integrals
can be evaluated using Theorem \ref{thmsym8}. For example:
\ba
\ift \frac{dz}{(z^{8} + 5z^{6} + 14z^{4} + 5z^{2} + 1)^{4} } & = &
\frac{(14325195794 + 2815367209 \, \sqrt{26} \, ) \, \pi}
{14623232 \, ( 9 + 2 \sqrt{26} \,)^{7/2} }. \nn
\ea

\vskip 5pt

\no
{\bf Example 8}. As in the case of degree $6$ we can provide numerical 
approximations to nonsymmetric integrals of degree $8$. The 
iteration (\ref{scheme6}) is now replaced by  \\
\ba
a_{1}^{n+1} & = & 
\frac{a_{2}^{n}(a_{1}^{n}+a_{3}^{n}) + 4a_{1}^{n}a_{3}^{n} + 
10(a_{1}^{n}+a_{3}^{n}) + 
8(a_{2}^{n}+2)}{(a_{1}^{n}+a_{2}^{n}+a_{3}^{n}+2)^{3/2} } \nn \\
a_{2}^{n+1} & = & \frac{a_{1}^{n}a_{3}^{n}+6(a_{1}^{n}+a_{3}^{n}) 
+ 2(a_{2}^{n}+10) }{a_{1}^{n}+a_{2}^{n}+a_{3}^{3}+2} \nn \\
a_{3}^{n+1} & = & \frac{a_{1}^{n}+a_{3}^{n}+8}{(a_{1}^{n} + a_{2}^{n} + a_{3}^{n}+2)^{1/2}} \nn \\
b_{0}^{n+1} & = & \frac{b_{0}^{n}+b_{1}^{n}+b_{2}^{n}+b_{3}^{n}}{(a_{1}^{n} + a_{2}^{n} + a_{3}^{n}+2)^{3/4}} \nn \\
b_{1}^{n+1} & = & \frac{b_{3}^{n}(3a_{1}^{n}+a_{2}^{n}+6)+b_{2}^{n}(a_{1}^{n}+4) +b_{1}^{n}(a_{3}^{n}+4) +b_{0}^{n}(3a_{3}^{n}+a_{2}^{n}+6) }
{(a_{1}^{n} + a_{2}^{n} + a_{3}^{n} +2)^{5/4} } 
\nn \\
b_{2}^{n+1} & = & \frac{b_{3}^{n}(a_{1}^{n}+5) + b_{2}^{n}+b_{1}^{n} + b_{0}^{n}(a_{3}^{n}+5) }{(a_{1}^{n} + a_{2}^{n} + a_{3}^{n}+2)^{3/4}} \nn \\
b_{3}^{n+1} & = & \frac{b_{0}^{n}+b_{3}^{n}}{(a_{1}^{n} + a_{2}^{n} + a_{3}^{n} +2)^{1/4}} \nn
\ea
\no
with initial conditions
$a_{1}^{0},a_{2}^{0},a_{3}^{0},b_{0}^{0},b_{1}^{0},b_{2}^{0},b_{3}^{0}$. 
Then 
\ba
U_{8}(a_{1},a_{2},a_{3},b_{0},b_{1},b_{2},b_{3}) & := & \int_{0}^{\infty} 
\frac{b_{0}x^{6} + b_{1}x^{4} + b_{2}x^{2} + b_{3}}{x^{8} + a_{1}x^{6} + 
a_{2}x^{4} + a_{3}x^{2} + 1} \; dx  \label{19u8}
\ea
\no
is invariant under these transformations.  \\

\no
{\bf Note}. Numerical calculations show that $(a_{1}^{n},a_{2}^{n},a_{3}^{n}) 
\to (4,6,4)$ and that $(b_{0}^{n},b_{1}^{n},b_{2}^{n},b_{3}^{n}) 
\to (1,3,3,1) \, L$ for some $L$ depending upon the initial conditions.

\vskip 10pt

\no
{\bf Example 9}. A symmetric function of degree $12$. We 
use Theorem \ref{main} to evaluate 
\ba
I & := & \ift \frac{z^{18} \, dz}{( z^{12} + 14z^{10} + 15z^{8} + 4z^{6} +
15z^{4} + 14z^{2} + 1)^{3}}  \nn
\ea
as
\ba
I & = & \frac{25 \pi ( 25 \sqrt{56} - 54) }{301989888}. \label{exam3}
\ea
\vskip 5pt

\no
Here $p =3, \; n=9$, and $m=2$, so $n > (m+1)p-1 $ and we need to 
apply the transformation $z \to 1/z$ to reduce the value of $n$. 
Indeed, we have 
\ba
I & = & \ift \frac{z^{16} \; dz}{ (z^{12} + 14z^{10} + 15z^{8} + 4z^{6} +
15z^{4} + 14z^{2} + 1)^{3}}, \nn
\ea
\no
and Theorem \ref{main} now yields
\ba
I & = & 2^{-17} \ift \frac{z^{16} \; dz}{ 131072 (1 + z^{2})^{3} 
( 1 + 4z^{2} + z^{4})^{3}}. \nn
\ea
\no
The new integrand is expanded in partial fractions in the variable $t = z^{2}$
to produce (\ref{exam3}). 

\vskip 15pt 

\no
{\bf Example 10}. We use Theorem \ref{main} to evaluate
\ba
I & := & \ift \frac{z^{10} \; dz }{Q^{2}(z)} \nn
\ea
\no
where
\ba
Q(z) & = & z^{20} + 6z^{18} + 93z^{16} -24z^{14} + 162z^{12}
+548z^{10} + 162z^8 -24z^6 +93z^4 +6z^2 +1. \nn
\ea
\no
The factorization 
\ba
Q(z) & = & ( 1+ z^{2})^2 T(z) T(-z)  \nn
\ea
\no
with 
\ba
T(z) & = & z^{8} -2z^{7} + 4z^{6} +14z^{5} +6z^{4} -14z^{3} + 4z^{2} + 2z +1
\nn
\ea
\no
leads to a partial fraction expansion containing the term
$$ 
\frac{72 -501 z + 1994z^2-2617z^3+1228z^4-43z^5 + 34z^6 -55z^7}
{8388608 ( 1 -2z+4z^2 +14z^3+6z^4 - 14z^5 + 4z^{6} + 2z^{7} +z^{8} )^{2} },   
$$
\no 
which we were unable to integrate; furthemore, the roots of $T(z)=0$ cannot 
be evaluated by radicals.  The 
procedure described in Theorem \ref{main}, however,  shows that
\ba
I & = & \ift \frac{z^{10} ( 4 + z^{2}) ( z^6 + 36z^4 +96z^2 + 64) \; dz }
{524288 (z^2 + 1)^2 \, (z^8 + 3z^6 + 8z^4 + 3z^2 +1)^2},  \nn
\ea 
\no
the integrand of which can be expanded to yield
\ba
I & = & - \frac{9}{8388608} \ift \frac{dz}{(z^{2} + 1)^{2} } 
- \frac{75}{8388608} \ift \frac{dz}{z^{2} + 1} \nn \\
 & & + \ift \frac{ 1921z^{6} + 10815z^4 +4111z^2 +1462} 
{2097152 \, (z^{8} + 3z^{6} + 8z^{4} + 3z^{2} + 1)^{2} } \; dz \nn \\
 & & \nn \\
 & & +  \ift \frac{ 91z^6 +719z^4 +1259z^2 -5764}
{8388608(z^{8} + 3z^{6} + 8z^{4} + 3z^{2} + 1)} \; dz. \nn
\ea
\no
Every piece is now computable, with the final result
$$ I =  \frac{(6480 - 509 \sqrt{15} ) \, \pi }{24159191040}.$$

\vskip 5pt 

\no
{\bf Example 11}. The symmetric functions of degree $16$ have denominator 
\ba
D_{4}(d_{1},d_{2},d_{3},d_{4};z) & = & z^{16} + d_{4}z^{14} + d_{3}z^{12} + 
d_{2}z^{10} +  2d_{1}z^{8} + d_{2}z^{6} + d_{3}z^{4} + d_{4}z^{2} + 1, \nn
\ea
the integral of which is computed in terms of 
\ba
E_{4}(d_{1},d_{2},d_{3},d_{4};z) & = & (1+d_{1}+d_{2}+d_{3}+d_{4})z^{8} + 
2(16+d_{2}+4d_{3}+ 9d_{4})z^{6}  \nn \\
 &  + & 8(20+d_{3}+6d_{4})z^{4} + 32(8+d_{4})z^{2} + 128. \nn
\ea
\no
This new integral is symmetric provided 
\ba
\begin{bmatrix}
d_{1} \\
d_{2} 
\end{bmatrix}
 & = & 
\begin{bmatrix}
15 \\
112 
\end{bmatrix}
+ 
\begin{bmatrix}
3 \\
-4 
\end{bmatrix}
d_{3} + 
\begin{bmatrix}
-8 \\
7 
\end{bmatrix}
d_{4}. 
\label{param16}
\ea
\no
Introduce the new parameters 
\ba
e_{j} = d_{j} - \binom{8}{5-j} \; \text{ for  } 2 \leq j \leq 4 \; 
\text{ and } 
e_{1}  =  d_{1} - \tfrac{1}{2} \binom{8}{4}. \nn
\ea
\no
Then (\ref{param16}) yields 
\ba
\begin{bmatrix}
e_{1} \\
e_{2} 
\end{bmatrix}
 & = & 
\begin{bmatrix}
3 \\
-4 
\end{bmatrix}
e_{3} + 
\begin{bmatrix}
-8 \\
7 
\end{bmatrix}
e_{4}. 
\label{newparam16}
\ea
\no
Thus, if the original denominator has the form 
\ba
D_{4}(z) & = & (z^{16} + 1 ) + d_{4}(z^{14}+z^{2}) + 
d_{3}(z^{12}+z^{4}) + (112 - 4d_{3}+7d_{4})(z^{10}+z^{6}) \nn \\
 &  & + 2(15+3d_{3}-8d_{4})z^{8}, \nn
\nn
\ea
\no
the integral
\ba
\ift \frac{P(z)}{\left( \; D_{4}(z) \; \right)^{m+1} } \, dz  \nn
\ea
\no
is reduced to an integral with symmetric denominator of degree $8$ and 
these are computable. We can thus evaluate 
a $2$-parameter family of symmetric  integrals 
of degree $16$.

For example, take $d_{3} = d_{4} = 1$ to obtain
\ba
R_{1}(z) & = &  \frac{z^{4}}{(z^{16} + z^{14} + z^{12} + 115z^{10} + 
20 z^{8} + 115z^{6} + z^{4} + z^{2} + 1)^{2}}. \label{fun1}
\ea
\no
The main theorem yields
\ba
R_{2}(z) & = & \frac{1024 z^{4} + 2304z^{6} + 1792 z^{8} + 560 z^{10} + 
60 z^{12} + z^{14} }{2^{7} ( 16z^{8} + 36z^{6} + 27z^{4} + 36z^{2} + 16)^{2}}
\label{fun2}
\ea
\no 
so that 
\ba
\ift R_{1}(z) \, dz & = & \ift R_{2}(z) \, dz. 
\nn
\ea
\no
Letting $f[n]:= N_{n,8}[1,n,27/32,9/4]$ we obtain
\ba
\ift R_{1}(z) \; dz & = & 2^{-15} \left( f[0]+60f[1]+1584f[2] + 
4096 f[3] \right) \nn
\ea
\no
and conclude that 
$$
\ift 
\frac{z^{4} \; dz }{(z^{16} + z^{14} + z^{12} + 115z^{10} + 
20 z^{8} + 115z^{6} + z^{4} + z^{2} + 1)^{2}}  \; =  $$
\vskip 2pt 
$$
= \; \frac{(149288517 + 12947003 \sqrt{131}) \pi}{1124663296 \sqrt{ 54925 + 
4798 \sqrt{131} } }. $$

\vskip 10pt

\no
{\bf Example 12}.  We classify the 
symmetric denominators of degree $32$ that yield computable integrals. 
These functions  depend on $8$ parameters
\ba
D_{8}(d_{1}, \cdots,d_{8};z) & = & \sum_{k=0}^{8} d_{9-k} 
(z^{2k} + z^{2(16-k)}) 
\ea
and the main theorem expresses the integral in terms of $E_{8}$. The 
conditions for $E_{8}$ to be 
symmetric yield 
\ba
d_{1} & = & -3441 + 35d_{5} +64d_{6} -312 d_{7} -3264 d_{8} \nn \\
d_{2} & = & 34720-56d_{5} -110d_{6} + 560d_{7} + 4565d_{8} \nn \\
d_{3} & = & -3472 + 28d_{5} + 64d_{6} -329d_{7} -2240d_{8} \nn \\
d_{4} & = & 4960 - 8d_{5} -19d_{6} + 80d_{7} + 938d_{8} \nn
\ea
\no
and the symmetric $E_{8}$ is 
\ba
E_{8}(d_{5},d_{6},d_{7},d_{8};z)  =  32768(1+z^{16}) + 
(131072 + 8192d_{8})(z^{2} + z^{14}) +   \nn
\ea
\ba
(212992 + 2048d_{7}+28672d_{8})(z^{4} + z^{12}) + 
(180224 + 512d_{6} + 6144d_{7} + 39424d_{8})(z^{6} + z^{10}) + \nn
\ea
\ba
+ (84480 + 128d_{5} + 1280d_{6} + 6912d_{7} + 26880d_{8})z^{8}. \nn
\ea
\no
The symmetry of $E_{8}$  now determines $d_{5},d_{6}$ in terms of 
$d_{7},d_{8}$ and we obtain
\ba
\begin{bmatrix}
d_{1} \\
d_{2} \\ 
d_{3} \\
d_{4} \\
d_{5} \\
d_{6} 
\end{bmatrix}
 & = & 
-31     
\begin{bmatrix}
63475 \\
-100800 \\ 
47936 \\
-13664 \\
2220 \\
- 224
\end{bmatrix}
 + 
\begin{bmatrix}
9166  \\
-14392\\ 
6895 \\
-1964 \\
322 \\
- 28
\end{bmatrix}
d_{7} + 
\begin{bmatrix}
54640  \\
-86645 \\ 
41664 \\
-11471  \\
2000  \\
-189
\end{bmatrix}
d_{8}. 
\label{parameter32}
\ea
\no
The function $(z^{2} + 1)^{16}$ is a symmetric polynomial of degree $32$ 
and yields a particular solution to (\ref{parameter32}). 

As before let 
\ba
e_{j}  =  d_{j} - \binom{16}{9-j} \; \text{ for  } 2 \leq j \leq 8 \; 
\text{ and }  
e_{1}  =  d_{1} - \tfrac{1}{2} \binom{16}{8}. \nn
\ea
\no
Then
\ba
\begin{bmatrix}
e_{1} \\
e_{2} \\ 
e_{3} \\
e_{4} \\
e_{5} \\
e_{6} 
\end{bmatrix}
 & = & 
\begin{bmatrix}
9166  \\
-14392\\ 
6895 \\
-1964 \\
322 \\
- 28
\end{bmatrix}
e_{7} + 
\begin{bmatrix}
54640  \\
-86645 \\ 
41664 \\
-11471  \\
2000  \\
-189
\end{bmatrix}
e_{8} 
\label{newparameter32}
\ea
\no
and as in the case of degree $16$ we can compute a $2$-parameter family of 
symmetric integrals of degree $32$. \\

\vskip 15pt

\no
\appendix
\section{Two binomial sums} \label{sec-binom}
\setcounter{equation}{0}

The closed-form evaluation of sums involving
binomials coefficients can be obtained by traditional analytical 
techniques or by using the powerful WZ-method as described in \cite{ab}.  We 
discuss two sums  used to simplify expressions in later sections, presenting 
one proof in each style. 

\begin{Lem}
\label{lem31}
Let $k, N$ be positive integers with $k \leq N$. Then 
\ba
\sum_{j=0}^{N-k} \binom{2N+1}{2j} \binom{N-j}{k} & = & 
\binom{2N-k}{k} 4^{N-k}. \label{identity1}
\ea
\end{Lem}

\begin{proof}
\no
Multiply the left hand side of (\ref{identity1}) by $x^{k}$ and sum over $k$ 
to produce
\ba
\sum_{k=0}^{N} \sum_{j=0}^{N-k} \binom{2N+1}{2j} \binom{N-j}{k} x^{k} 
& = & 
\sum_{j=0}^{N} \binom{2N+1}{2j} \sum_{k=0}^{N-j} \binom{N-j}{k} x^{k}  \nn \\
& = & \sum_{j=0}^{N} \binom{2N+1}{2j} (x+1)^{N-j} \nn \\
& = & \frac{( 1 + \sqrt{x+1})^{2N+1} 
- (1 - \sqrt{x+1} )^{2N+1} }{2 \sqrt{x+1} } \nn \\
& = & \sum_{k=0}^{N} \binom{2N-k}{k} 4^{N-k} x^{k}. \nn
\ea
\no
In order to justify the last step we start with the well known result
\ba
\frac{1}{\sqrt{1 - 4y}} \left( \frac{1 - \sqrt{1-4y}}{2y} \right)^{i} 
& = & \sum_{k} \binom{2k+ i}{k} y^{k} \label{wellknown}
\ea
\no
(see WILF \cite{wilf}, page 54).  Letting $x = -4y$ and $i = 2N+1$ 
in (\ref{wellknown}) gives 
\ba
\frac{(1- \sqrt{x+1} \,)^{2N+1} }{2 \sqrt{x+1}} & = & 
\sum_{k=0}^{\infty} (-1)^{k+1} \binom{2N+1+2k}{k} 4^{-(N+1+k)} x^{2N+1+k};
\nn
\ea
\no
similarly $x = -4y$ and $i = -2N-1$  yields
\ba
\frac{(1+ \sqrt{x+1} \,)^{2N+1} }{2 \sqrt{x+1}} & = & 
\sum_{k=2N+1}^{\infty} (-1)^{k+1} \binom{2N+1+2k}{k} 4^{-(N+1+k)} x^{2N+1+k}. 
\nn
\ea
\no
Thus 
\ba
\frac{( 1 + \sqrt{x+1} \,)^{2N+1} 
- (1 - \sqrt{x+1}  \,)^{2N+1} }{2 \sqrt{x+1} }  & = & 
\sum_{k=0}^{2N} (-1)^{k} \binom{-(2N+1-2k)}{k} 4^{N-k} x^{k}  \nn \\
& = & \sum_{k=0}^{N} \binom{2N-k}{k} 4^{N-k} x^{k}.  \nn 
\ea
\end{proof}

\vskip 5pt 

\no
\begin{Lem}
\label{lem32}
Let $N \in \mathbb{N}$. Then 
\ba
\sum_{j=0}^{N} \binom{2N+1}{2j} (1+z^{2})^{N-j} 
& = & \sum_{j=0}^{N} \binom{N+j}{2j} 4^{j}  z^{2(N-j)}. \nn
\ea
\end{Lem}

\no
\begin{proof}
The coefficient of $z^{2k}$ on the left hand side is 
\ba
\sum_{j=0}^{N-k} \binom{2N+1}{2j} \binom{N-j}{k}, \nn
\ea
\no
and the corresponding coefficient on the right hand side is
$\binom{2N-k}{k} 4^{N-k}$. The result then follows from Lemma \ref{lem31}. 
\end{proof}

\begin{Lem}
\label{lem33}
Let $k,N \in \mathbb{N}$ with $k \leq N$. Then
\ba
\sum_{j=0}^{k} \binom{2N}{2j} \binom{N-j}{N-k} & = & \frac{2^{2k-1} N}{k} 
\binom{k+N-1}{N-k}. \nn
\ea
\end{Lem}

\begin{proof}
This lemma could be proven in the same style as Lemma \ref{lem31}. Instead we
use the WZ-method as explained in \cite{ab}. 
Indeed, let
\ba
F(k;j) & = & \frac{k \binom{2N}{2j} \binom{N-j}{N-k} }
{ N 2^{2k-1} \binom{k+N-1}{N-k} },  \nn
\ea
\no
and define, with the package EKHAD, the function
\ba
G(k;j) & = & F(k;j) \times \frac{j(2j-1)}{2(N+k)(k-j+1)}. \nn
\ea
\no
Then $F(k;j) - F(k+1;j)  =  G(k;j+1) - G(k;j)$, and 
summing over $j$ we see that the sum of $F(k;j)$ over $j$ 
is independent of $k$.  The case $k=N$ produces $1$ as the common value.
\end{proof}
\medskip

\no
\begin{Lem}
\label{lem34}
Let $p \in {\mathbb{N}}, \, d_{1}, d_{2}, \cdots, d_{p}$ be 
parameters, and define
$d_{p+1} := 1$. Then 
\ba
\sum_{k=0}^{p} d_{p+1-k} z^{2k} \; \sum_{j=0}^{p-k}  
\binom{2p-2k}{2j} (1 + z^{2})^{p-k-j}  \; \; = \label{star}
\ea
\ba
\left( \sum_{j=1}^{p+1} d_{j} \right) z^{2p} \; + \; \; 
\sum_{i=1}^{p} 2^{2i-1} z^{2(p-i)} \left( \sum_{j=1}^{p+1-i} 
\frac{j+i-1}{i} \; \binom{j+2i-2}{j-1} \, d_{j+i} \right). \nn
\ea
\end{Lem}

\begin{proof}
For fixed $0 \leq i \leq p-1$  the coefficient of $z^{2 i}$ on the 
right hand side of (\ref{star}) is
\ba
\left[ {\rm{RHS}} \right] (2 i) & = & \frac{2^{2(p- i)-1}}{p- i} 
\sum_{r=p- i +1}^{p+1} ( r-1) \binom{r+p- i -2}{r - p + i -1} \; d_{r}, 
\nn
\ea
\no
and for $ i = p$ we have
$\left[ {\rm{RHS}} \right](2p)  =  1 + \sum_{j=1}^{p} d_{j}$. 
Similarly, for the left hand side of (\ref{star}),
\ba
\left[ {\rm{LHS}} \right] (2 i ) & = & 
\sum_{r=p+1- i}^{p+1} d_{r} 
\left( \sum_{j=0}^{p- i} \binom{2r-2}{2j} \binom{r-j-1}{r-j-1+ 
i} \right). \nn
\ea
\no
It is easy to check that the coefficients 
of $z^{2p}$ match. It suffices to 
show that for each 
$i$ such that $ 0 \leq i \leq p-1$ and for each $r$ such that 
$ p+1 - i \leq r \leq p+1$ we have
\ba
\sum_{j=0}^{p - i} \binom{2r-2}{2j} \binom{r-1-j}{r-1-p + i}  & = &
\frac{2^{2(p- i) - 1}}{p- i}
(r -1) \binom{r + p - i - 2}{r - p + i-1}. \nn
\ea
\no
This follows from Lemma \ref{lem31} with $k = p - i$ and $N = r - 1$.
\end{proof}

\thanks{The suggestions of the referrees and the editor are gratefully acknowledged.}


\begin{thebibliography}{99}

\bibitem{unim}
BOROS, G. - MOLL, V.: {\em A criterion for unimodality}. Elec. Jour. of 
Combinatorics, {\bf 6} (1999), \#R10.  

\bibitem{hyp}
BOROS, G. - MOLL, V.: {\em An integral hidden in Gradshteyn and Ryzhik}. 
Jour. Comp. Appl. Math. {\bf 106}, 361-368, 1999.

\bibitem{landen}
BOROS, G. - MOLL, V.: {\em A rational Landen transformation}. 
Contemporary Mathematics {\bf 251}, 83-91, 2000.

\bibitem{borbor}
BORWEIN, J. - BORWEIN, P.: {\em Pi and the AGM}.  Canadian Mathematical 
Society. Wiley-Interscience Publication. 

\bibitem{bron}
BRONSTEIN, M.: {\em Symbolic Integration I. Transcendental functions}. 
Algorithms and Computation in Mathematics, {\bf 1}. Springer-Verlag, 1997.

\bibitem{geddes}
GEDDES, K. - CZAPOR, S.R. - LABAHN, G.: {\em Algorithms for Computer Algebra}.
Kluwer, Dordrecht. The Netherlands, 1992.

\bibitem{gauss1}
GAUSS, K.F.: {\em Arithmetische Geometrisches Mittel}, 1799. In {
\em Werke}, {\bf 3}, 
361-432. Konigliche Gesellschaft der Wissenschaft, Gottingen. Reprinted 
by Olms, Hildescheim, 1981.

\bibitem{gr}
GRADSHTEYN, I.S. - RYZHIK, I.M.: {\em Table of Integrals, Series and Products}.
Fifth Edition, ed. Alan Jeffrey. Academic Press, 1994.

\bibitem{hardy}
HARDY, G.H.: {\em The Integration of Functions of a Single Variable}. Cambridge
Tracts in Mathematics and Mathematical Physics, 2, Second Edition, Cambridge 
University Press, 1958.

\bibitem{herm}
HERMITE, C.: {\em Sur l'integration des fractions rationelles}. Nouvelles 
Annales de Mathematiques ( $2^{\text{eme}}$ serie) {\bf 11}, 145-148, 1872.

\bibitem{horo}
HOROWITZ, E.: {\em Algorithms for partial fraction decomposition and rational 
function integration}. Proc. of SYMSAM'71, ACM Press, 441-457, 1971.

\bibitem{lazrio}
LAZARD, D. - RIOBOO, R.: {\em Integration of Rational Functions: Rational 
Computation of the Logarithmic Part}. Journal of Symbolic Computation {\bf 9}, 
113-116, 1990.

\bibitem{mckmoll}
MCKEAN, H. - MOLL, V.: {\em Elliptic Curves:  Function Theory, Geometry, 
Arithmetic}. Cambridge University Press, 1997.

\bibitem{newman}
NEWMAN, D.: {\em A simplified version of the fast algorithm of Brent and
Salamin}. Math. Comp. {\bf 44}, 207-210, 1985.

\bibitem{ostro}
OSTROGRADSKY, M.W.: {\em De l'integration des fractions rationelles}. Bulletin
de la Classe Physico-Mathematiques de l'Academie Imperieriale des Sciences
de St. Petersbourgh, IV, 145-167, 286-300. 1845.

\bibitem{ab}
PETKOVSEK, M. - WILF, H.S. - ZEILBERGER, D.: {\em A=B}. A. K. Peters, 
Wellesley, Massachusetts. 1996.

\bibitem{roth}
ROTHSTEIN, M.: {\em A new algorithm for the integration of Exponential and 
Logarithmic Functions}, Proc. of the 1977 MACSYMA Users Conference, NASA
Pub., CP-2012, 263-274.

\bibitem{trag}
TRAGER, B.M.: {\em Algebraic factoring and rational function integration}. 
Proc. SYMSAC 76, 219-226.

\bibitem{wilf}
WILF, H.S.: {\em generatingfunctionology}. Academic Press, 1990.

\end{thebibliography}
\end{document}